\documentclass[11pt]{amsart}

\theoremstyle{plain}
\newtheorem{theorem}{Theorem}
\newtheorem{definition}{Definition}
\newtheorem{lemma}{Lemma}
\newtheorem{corollary}{Corollary}
\newtheorem{remark}{Remark}
\newtheorem{example}{Example}

\begin{document}

\title[On The Development of Nonlinear Operator Theory]
{On The Development of Nonlinear Operator Theory}
\author{Wen Hsiang Wei}
\address{Department of Statistics, Tung Hai University, Taiwan}
\email{wenwei@thu.edu.tw}
\subjclass[2010]{Primary 47H99; Secondary 46H30}

\begin{abstract}

The basic results for nonlinear operators are given. These results include
nonlinear versions of classical uniform boundedness theorem and
Hahn-Banach theorem. Furthermore, the mappings from a
metrizable space into another normed space can fall in some normed spaces by defining
suitable norms. The results for the mappings on the
metrizable spaces can be applied to the operators on the space of bounded linear
functionals corresponding to the Dirac's delta function.

\end{abstract}

\keywords{Banach algebra, Dirac's delta function,
Nonlinear Hahn-Banach theorem, Nonlinear operators,
Nonlinear uniform boundedness theorem}

\maketitle

\section{Introduction}

Operator theory has been at the heart of research in analysis (see
\cite{abramovich}; \cite{pier}, Chapter 4). Moreover, as implied by
\cite{neuberger}, considering nonlinear case should be
essential. Developing useful results for the operators holds the promise for the
wide applications of nonlinear functional analysis to a variety of scientific areas.

In classical functional analysis, the space of bounded linear operators is a
normed space endowed with a sensible norm.
In next section, several classes of normed functions are defined and
some set including certain nonlinear operators from a normed space into a normed space
turns out to be a normed space.

For the bounded linear operators between the normed spaces, the
Hahn-Banach theorem and the uniform boundedness theorem are basic
theorems. It is sensible to develop the nonlinear counterparts of these
theorems. The nonlinear uniform boundedness theorem and the nonlinear
Hahn-Banach
theorem are given in Section 3 and Section 4, respectively.

Some mappings, for example, the distributions on some
metrizable spaces, do not have the "normed" values because the metrizable
spaces are not normable. To resolve the problem, the mappings from a
metrizable space into a normed space can have the normed values by defining a
suitable norm depending on the metric of the metrizable space. The
results for the mappings on the metrizable spaces along with
some examples are given in the last section.

\section{Nonlinear functional spaces}

Let \( X \) and \( Y \) be the normed
spaces over the field \( K \) with some sensible norms \( ||\cdot||_{X} \) and  \( ||\cdot||_{Y} \),
respectively, where \( K \) is either a real field \( R \) or a complex field \( C \).
Note that the vector spaces and the normed spaces in this
article are assumed to be not trivial, i.e., not only including the zero element.
Let
\( V(X,Y) \) be the set of all operators from \( X  \) into \( Y \), i.e.,
the set of arbitrary maps from \( X  \) into \( Y \). Note that the operators
in \( V(X,Y) \) are not assumed to be continuous.
Let the algebraic operations of \( F_{1}, F_{2} \in V(X,Y) \) be the operators
from \( X \) into
\( Y \) with \( (F_{1}+F_{2})(x)=F_{1}(x)+F_{2}(x) \) and
\( (\alpha F_{1})(x)=\alpha F_{1}(x) \) for \( x \in X \), where \( \alpha \in K \) is a
scalar. Also let the
zero element in
\( V(X,Y) \) be the operator with the image equal to the zero element in \( Y \).
\( V(X,Y) \) is a vector space over the field \( K \). Define a
non-negative extended real-valued function
\( p \), i.e., the range of \( p \) including \( \infty \), on \( V(X,Y) \) by
\begin{eqnarray*}
p(F)=\max\left(\sup_{x \neq 0,x \in X} \frac{\|F(x)\|_{Y}}{\|x\|_{X}},
\|F(0)\|_{Y}\right)
\end{eqnarray*}
for \( F \in V(X,Y) \). The non-negative extended real-valued function \( p \) is a generalization of the
norm for the linear operators. Let
\( B(X,Y) \), the subset of \( V(X,Y) \), consist of all operators with
\( p(F) \) being finite. Note that
\( p \) is a norm on \( B(X,Y) \) and \( [B(X,Y),p] \) is a normed space.

\begin{remark}

For \( F \in B(X,Y) \), another norm equivalent to \( p \) is
\begin{eqnarray*}
p^{\ast}(F)=\sup_{x \neq 0,x \in X} \frac{\|F(x)\|_{Y}}{\|x\|_{X}} +\|F(0)\|_{Y}.
\end{eqnarray*}
In addition, two classes of norms related to \( p \) and \( p^{\ast} \) are
\begin{eqnarray*}
q(F)=\max\left(\sup_{x \neq 0,x \in X} \frac{\|F(x)\|_{Y}}{\|x\|_{X}},
s\|F(0)\|_{Y}\right)
\end{eqnarray*}
and
\begin{eqnarray*}
q^{\ast}(F)=\sup_{x \neq 0,x \in X} \frac{\|F(x)\|_{Y}}{\|x\|_{X}} +s\|F(0)\|_{Y},
\end{eqnarray*}
\( s >  0 \). Furthermore, the other class of the
non-negative extended real-valued functions on \( V(X,Y) \) which
includes the function \( p \) is
\begin{eqnarray*}
p_{k}(F)=\max\left(\sup_{x \neq 0, x \in X}
\frac{\left|\left|F(x)\right|\right|_{Y}}{\left|\left|x\right|\right|_{X}^{k}},
\left|\left|F(0)\right|\right|_{Y}\right),
\end{eqnarray*}
where \( k\) is a positive number. Note that
\begin{eqnarray*}
\left|\left|F(x)\right|\right|_{Y} \leq p_{k}(F) \left|\left|x\right|\right|^{k}_{X},
x \neq 0,
\end{eqnarray*}
i.e., the normed value of \( F(x) \) being dominated by the power of the normed
value of \( x \neq 0 \), if \( p_{k}(F) \) is finite.

In addition, let \( V(S,Y) \) be the set of all operators from the set \( S \subset X  \) into \( Y \) and
\( 0 \in S \). Also let the zero element in
\( V(S,Y) \) be the operator of which image equal to the zero element in \( Y \). Then
\( V(S,Y) \) is a vector space.  Let \( B(S,Y) \) be the subset of \( V(S,Y) \)
with the property that
\( ||F||_{B(S,Y)} \) is
finite for all \( F \in B(S,Y) \), where
\begin{eqnarray*}
\left|\left|F\right|\right|_{B(S,Y)}=\max\left(\sup_{x \neq 0,x \in S}
\frac{\|F(x)\|_{Y}}{\|x\|_{X}},
\|F(0)\|_{Y}\right).
\end{eqnarray*}
Note that \( B(S,Y) \) is a normed space, i.e., \( ||\cdot||_{B(S,Y)} \) being a
normed function on
\( B(S,Y) \).

\end{remark}

Hereafter the norm \( p \) is used, i.e., \( ||F||_{B(X,Y)}=p(F) \) for \( F \in
B(X,Y)
\). Note that the bounded linear operators fall in
\( B(X,Y) \). However, unlike a linear operator,
a continuous nonlinear operator might not fall in the space \( B(X,Y) \).
As \( X=Y \), the notation
\( B(X)=B(X,X) \) is used. If the other norm is used, for example, \( p_{k} \),
the notation \( [B(X,Y),p_{k}] \) specifying the norm is used.

Let the notation of the composition of two operators
be \( \circ \) in the following lemma.

\begin{lemma}

Let \( F_{1} \in B(X,Y) \) and \( F_{2} \in B(Y,Z) \), where \( X , Y, \) and \( Z \)  are normed
spaces. Then
\begin{eqnarray*}
\|F_{1}(x)\|_{Y} \leq \|F_{1}\|_{B(X,Y)}\|x\|_{X},x \neq 0.
\end{eqnarray*}
 If \( F_{1}(0)=F_{2}(0)=0 \), then
\begin{eqnarray*}
\|F_{2} \circ F_{1}\|_{B(X,Z)} \leq \|F_{1}\|_{B(X,Y)}\|F_{2}\|_{B(Y,Z)},
\end{eqnarray*}
i.e., \( F_{2} \circ F_{1} \in B(X,Z) \).

\end{lemma}

\textsc{Proof.} As \( x \neq 0 \),
\begin{eqnarray*}
\frac{\|F_{1}(x)\|_{Y}}{\|x\|_{X}} \leq
\sup_{x \neq 0, x \in X} \frac{\|F_{1}(x)\|_{Y}}{\|x\|_{X}} \leq \|F_{1}\|_{B(X,Y)}
\end{eqnarray*}
and hence \( ||F_{1}(x)||_{Y} \leq ||F_{1}||_{B(X,Y)}||x||_{X} \). In addition, as
\( F_{1}(0)=F_{2}(0)=0 \),
\begin{eqnarray*}
\|F_{2} \circ F_{1}\|_{B(X,Z)}
&=& \sup_{x \neq 0, x \in X} \frac{\|(F_{2} \circ F_{1})(x)\|_{Z}}{\|x\|_{X}} \\
&\leq& \sup_{x \neq 0, x \in X}
\frac{\|F_{2}\|_{B(Y,Z)}\|F_{1}(x)\|_{Y}}{\|x\|_{X}} \\ &=&
\|F_{2}\|_{B(Y,Z)}\left(\sup_{x \neq 0, x \in X} \frac{\|F_{1}(x)\|_{Y}} {\|x\|_{X}}\right)
\\ &=& \|F_{2}\|_{B(Y,Z)} \|F_{1}\|_{B(X,Y)}.
\end{eqnarray*}

\begin{flushright}
$\diamondsuit$
\end{flushright}

\section{Nonlinear uniform boundedness theorem}

For a family of bounded linear operators, the pointwise boundedness implies
the uniform boundedness. Theorem 3, the main theorem in this section and
considered as the nonlinear version of the uniform boundedness theorem,
gives the analogous result for certain nonlinear operators.

The bounded linear operators or the bounded nonlinear operators between normed spaces map
bounded sets into bounded sets. The class of operators mapping the bounded
sets into the bounded sets is defined as follows.
\begin{definition}

An operator \( F: X \rightarrow Y \) is topology bounded if and only if \( F \)
maps bounded sets in the normed space \( X \) into bounded sets in the
normed space \( Y
\).

\end{definition}

The operator \( F \in B(X,Y) \) is also called the norm bounded operator. A
linear operator is norm bounded if and only if it is topology bounded. It might
not be true for the nonlinear operator. A nonlinear operator is topology
bounded might not be norm bounded. However, given some sufficient
condition, a topology bounded nonlinear operator can be norm bounded, as
indicated by the following theorem.

\begin{theorem}

If \( F \in B(X,Y) \), then \( F \) is topology bounded. On the other hand, if \( F \)
is topology bounded and
\begin{eqnarray*}
\left|\left|F(kx)\right|\right|_{Y} \leq M\left|k\right|\left|\left|F(x)\right|\right|_{Y}
\end{eqnarray*}
for some positive constant \( M \), any scalar \( k \neq 0 \), and any \( x \neq 0\), then \( F \in B(X,Y)
\).

\end{theorem}

\textsc{Proof.} If \( F \) is norm bounded, then for any bounded set
\( V \in X \), there exists an open ball \( B_{r} \in X \) with a radius \( r > 0 \) and
a center at \( 0 \) such that \( V \subset B_{r} \) and \( \sup_{x \in
B_{r}}||F(x)||_{Y}
\leq ||F||_{B(X,Y)}\max(r,1) \) by Lemma 1, i.e., \( F(V) \) being bounded. Conversely,
let \( \partial \bar{B}_{1} \) be the
boundary of the unit closed ball centered at \( 0 \) in \( X \). Because \( F \) is
topology bounded, there exists some nonnegative number \(
\overline{k} \) such that
\(
\overline{k}=sup_{x
\in
\partial
\bar{B}_{1}} ||F(x)||_{Y} \). Then for any \( x \in X, x \neq 0 \),
\begin{eqnarray*}
\frac{\left|\left|F(x)\right|\right|_{Y}}{\left|\left|x\right|\right|_{X}}
\leq \frac{M\left|\left|x\right|\right|_{X}
\left|\left|F(x/\left|\left|x\right|\right|_{X})\right|\right|_{Y}}
{\left|\left|x\right|\right|_{X}}
\leq M\overline{k}.
\end{eqnarray*}
Therefore, \( ||F||_{B(X,Y)} \leq \max(M\overline{k},||F(0)||_{Y}) \) and \( F
\in B(X,Y) \).

\begin{flushright}
$\diamondsuit$
\end{flushright}

The inequality in the above theorem, referred to as the $M$-contraction property,
turns out to be crucial for the uniform boundedness of a certain family of
nonlinear operators.

\begin{definition}

An operator \( F \) from the normed space \( X \) into the normed space \( Y \)
is called a $M$-contraction operator if and only if
\begin{eqnarray*}
\left|\left|F(kx)\right|\right|_{Y} \leq M\left|k\right|\left|\left|F(x)\right|\right|_{Y}
\end{eqnarray*}
for some positive constant \( M  \), any scalar \( k \neq 0 \), and any \( x \neq 0 \). The set of all
topology bounded $M$-contraction operators is denoted as
\( B^{M}(X,Y)
\).

\end{definition}

By Theorem 1, a $M$-contraction operator is topology bounded if and only if it
is norm bounded.

\begin{theorem}

\( B^{M}(X,Y) \) is a closed subset of \( B(X,Y) \).

\end{theorem}

\textsc{Proof.}  By Theorem 1, a topology bounded $M$-contraction operator is norm
bounded, i.e., \( B^{M}(X,Y) \subset B(X,Y) \). Suppose
\( F_{n} \mathop{\longrightarrow}\limits_ {n \rightarrow\infty} F \), where
\( F_{n} \in B^{M}(X,Y) \) and \( F \in B(X,Y) \). Then there
exists \( m \) such that
\begin{eqnarray*}
&& \left|\left|F(kx)\right|\right|_{Y} \\ &\leq&
\left|\left|F_{m}(kx)-F(kx)\right|\right|_{Y}+
\left|\left|F_{m}(kx)\right|\right|_{Y} \\
&\leq& \epsilon +M\left|k\right|\left(
\left|\left|F_{m}(x)-F(x)\right|\right|_{Y}+
\left|\left|F(x)\right|\right|_{Y}\right) \\
&\leq&(M\left|k\right|+1)\epsilon +M\left|k\right|\left|\left|F(x)\right|\right|_{Y}
\end{eqnarray*}
for every scalar \( k \neq 0 \), \( x \neq 0 \) and any \( \epsilon > 0 \). Therefore,
\( F \in B^{M}(X,Y) \) and \( B^{M}(X,Y) \) is closed.

\begin{flushright}
$\diamondsuit$
\end{flushright}

Note that the space of the bounded linear operators is the subset of
\( B^{M}(X,Y) \) for any \( M \geq 1 \). The uniform boundedness property of the nonlinear operators
of interest is described as follows.

\begin{definition}

Let \( \{F_{\alpha}\} \), a subset of \( V(X,Y) \), be a family of operators, where \( \alpha \in I
\) and
\( I
\) is an index set. \(
\{F_{\alpha}\} \) is uniformly topology bounded if and only if
for any bounded set \( U \in X \) there
exists a bounded set \( V \in Y \) satisfying \( F_{\alpha}(U) \subset V \) for
every
\(
\alpha
\in I \). \( \{F_{\alpha}\} \) is uniformly norm bounded if and only if \( ||F_{\alpha}||_{B(X,Y)}
\leq c \) for every \(
\alpha
\in I \) and some positive constant \( c \).

\end{definition}

The uniform boundedness theorem holds for certain nonlinear operators which
have the $M$-contraction property.

\begin{theorem}

Let \( \{F_{\alpha}\} \) be a family of continuous $M$-contraction operators
from a Banach space \( X \) into a normed space
\( Y \), where \( \alpha \in I \) and \( I \) is an index set. \( \{F_{\alpha}\} \)
is uniformly norm bounded if the following conditions hold. \\
 (a) \( \{||F_{\alpha}(x)||_{Y}:\alpha \in I\} \) is bounded for \( x \in X \), i.e.,
\( ||F_{\alpha}(x)||_{Y} \leq c_{x} \) for every \( \alpha \in I \), where
\( c_{x} \) depending on \( x \)  is a positive number.\\
(b) There exists a positive constant \( L \) such that
\begin{eqnarray*}
\left|\left|F_{\alpha}(x_{1}+x_{2})\right|\right|_{Y} \leq
L\left|\left|F_{\alpha}(x_{1})+F_{\alpha}(x_{2})\right|\right|_{Y}
\end{eqnarray*}
for every \( \alpha \in I \) and \( x_{1},x_{2} \in X \).

\end{theorem}

\textsc{Proof.} Let \( A_{k}=\{x:||F_{\alpha}(x)||_{Y} \leq k,x \in X,\alpha \in I\}, k=1,2,\ldots. \).
\( A_{k} \) is closed by the continuity of \( F_{\alpha} \) and
\( X=\cup_{k=1}^{\infty}A_{k} \). Since \( X \) is of second category, there
exists \( A_{k_{0}} \) such that \( B(x_{0},r) \subset A_{k_{0}} \), where \( k_{0} \) is some positive integer and
\(  B(x_{0},r) \) is an open ball with the center \( x_{0} \in X \) and the radius \( r > 0 \). For
any \( x \neq 0 \), there exists a constant \( c > 1 \) such that \(
x=cr^{-1}||x||_{X}(z-x_{0})
\), where \( z \in B(x_{0},r) \) and \( z \neq x_{0} \). Then
\begin{eqnarray*}
&& \left|\left|F_{\alpha}(x)\right|\right|_{Y} \\ &\leq&
Mcr^{-1}||x||_{X}\left|\left|F_{\alpha}(z-x_{0})\right|\right|_{Y} \\ &\leq&
MLcr^{-1}||x||_{X}\left|\left|F_{\alpha}(z)+F_{\alpha}(-x_{0})\right|\right|_{Y}
\\ &\leq& MLcr^{-1}||x||_{X}\left(
k_{0}+M\left|\left|F_{\alpha}(x_{0})\right|\right|_{Y}\right)\\ &\leq&
M(M+1)Lcr^{-1}k_{0}||x||_{X}.
\end{eqnarray*}
Therefore,
\begin{eqnarray*}
&& \left|\left|F_{\alpha}\right|\right|_{B(X,Y)} \\ &\leq&
\max \left[M(M+1)Lcr^{-1}k_{0},c_{0} \right]
\end{eqnarray*}
for every \( \alpha \in I \) and \( \{F_{\alpha}\} \) is uniformly norm bounded.

\begin{flushright}
$\diamondsuit$
\end{flushright}

A bounded linear operator is a $M$-contraction operator and condition (b) in
the above theorem also holds for the bounded linear operator, i.e., Theorem 3
being a generalization of the classical
uniform boundedness theorem.

\begin{corollary}

Let \( F: X \rightarrow Y \), \( \{F_{n}\} \) be a sequence of continuous
$M$-contraction operators from a Banach space \( X \) into a normed space
\( Y \) and \( \{F_{n}(x)\} \) converges to \( F(x) \) with respect to the norm topology \(
||\cdot||_{Y} \) for every \( x \in X \). If there exists a positive constant \( L \) such that
\begin{eqnarray*}
\left|\left|F_{ n}(x_{1}+x_{2})\right|\right|_{Y} \leq
L\left|\left|F_{n}(x_{1})+F_{n}(x_{2})\right|\right|_{Y}
\end{eqnarray*}
for every \( n \) and \( x_{1},x_{2} \in X \), then the sequence
\( \{F_{n}\} \) is uniformly norm bounded and
\( F \in B(X,Y) \).

\end{corollary}

\textsc{Proof.} Because
\begin{eqnarray*}
\left|\left|F_{n}(x)\right|\right|_{Y} \leq \left|\left|F_{n}(x)-F(x)\right|\right|_{Y}
+\left|\left|F(x)\right|\right|_{Y}
\end{eqnarray*}
and \( ||F_{n}(x)-F(x)||_{Y} \mathop{\longrightarrow}\limits_ {n
\rightarrow\infty}  0 \) for every \( x \in X \), hence \( \{F_{n}(x)\} \) is bounded for every \( x \in
X \). By Theorem 3, \(
\{F_{n}\} \) is uniformly norm bounded, \i.e., \( ||F_{n}||_{B(X,Y)} \leq c \) for
every \( n \) and some positive constant \( c  \). Finally, for \( x \neq 0 \),
\begin{eqnarray*}
\left|\left|F(x)\right|\right|_{Y} =\lim_{n \rightarrow \infty}\left|\left|F_{n}(x)\right|\right|
\leq c \left|\left|x\right|\right|_{X}
\end{eqnarray*}
by Lemma 1 and
\begin{eqnarray*}
\left|\left|F(0)\right|\right|_{Y} =\lim_{n \rightarrow \infty}\left|\left|F_{n}(0)\right|\right|
\leq c,
\end{eqnarray*}
hence \( ||F||_{B(X,Y)} \leq c \), i.e.,
\( F \in B(X,Y) \).

\begin{flushright}
$\diamondsuit$
\end{flushright}

\section{Nonlinear Hahn-Banach theorems}

The Hahn-Banach theorem
states that a linear functional on a subspace of a vector space \( X \) can be extended
to the whole space with two properties preserved, linearity and the inequality
for the linear functional and a sub-linear functional. It turns out that a nonlinear
functional on a subset of a vector space can be extended to the whole space
with a certain inequality preserved, as given in Theorem 4. Further, as \( X \) is
a separable Hilbert space, Theorem 5 states that the extension to the whole
space also holds with both the inequality similar to the one in Theorem 4 and
the continuity of the nonlinear functional preserved. Theorem 6, the last
theorem in this section, is concerned with the extension results for the
nonlinear functionals with some specific forms. Hereafter, let
\( Dom(F) \) be the domain of  \( F \) which is the operator or the functional.

\begin{theorem}

Let \( X \) be a vector space and \( p \) be a sub-additive functional on
\( X
\). Let
\( F: S \rightarrow R \) be a functional on \( S \) and
\begin{eqnarray*}
F(s_{1})+F(s_{2}) \leq p(s_{1}+s_{2}),s_{1},s_{2} \in S, s_{1}
\neq s_{2},
\end{eqnarray*}
where \( S \) is a proper subset of \( X \).
\( F \) has an extension \( \hat{F}: X \rightarrow R \) satisfying
\begin{eqnarray*}
\hat{F}(s)=F(s), s \in S,
\end{eqnarray*}
and
\begin{eqnarray*}
\hat{F}(x_{1})+\hat{F}(x_{2}) \leq p(x_{1}+x_{2}),x_{1},x_{2} \in X,
x_{1} \neq x_{2}.
\end{eqnarray*}

\end{theorem}

\textsc{Proof.} By Zorn's lemma, the set of all extensions of \( F \) satisfying
the inequality has a maximal element by defining the partial ordering as the
inclusion of the domains of the extensions. The maximal element \( \hat{F} \)
being defined on the whole space \( X \) is proved next.

Suppose that \( Dom(\hat{F}) \) is a proper subset of \( X \). Since
\begin{eqnarray*}
\hat{F}(x_{1})+\hat{F}(x_{2}) \leq p(x_{2}+x_{2}) \leq
p(x_{1}-y)+p(x_{2}+y),
\end{eqnarray*}
\begin{eqnarray*}
\hat{F}(x_{2})-p(x_{2}+y) \leq p(x_{1}-y)-\hat{F}(x_{1})
\end{eqnarray*}
thus, where \( y \in X/Dom(\hat{F}) \), i.e., \( y \) falling in the intersection
of \( X \) and the complement of \( Dom(\hat{F}) \), and \( x_{1},x_{2} \in Dom(\hat{F}), x_{1} \neq x_{2}
\). Define \( F_{1}: Dom(\hat{F}) \cup \{y\} \rightarrow R \) by \( F_{1}(x)=\hat{F}(x)
\) for \( x \in Dom(\hat{F}) \) and
\( F_{1}(y)=-c \), where
\begin{eqnarray*}
c=\sup_{x \in Dom(\hat{F})} \left[\hat{F}(x)-p(x+y)\right]
\end{eqnarray*}
and the existence of \( c \), i.e., the supremum being finite,  is due to
\begin{eqnarray*}
c \leq \max\left[p(x_{1}-y)-\hat{F}(x_{1}),p(x_{2}-y)-\hat{F}(x_{2})\right].
\end{eqnarray*}
Then for \( x \in Dom(\hat{F}) \),
\begin{eqnarray*}
F_{1}(x)+F_{1}(y)=\hat{F}(x)-c \leq p(x+y).
\end{eqnarray*}
Thus, \( F_{1} \) satisfying the inequality is an extension of \( \hat{F} \), i.e., a
contradiction.

\begin{flushright}
$\diamondsuit$
\end{flushright}

The Hahn-Banach extension theorem for the linear functionals is a special case
of the following corollary.

\begin{corollary}

Let \( X \) be a vector space, \( S \) be a proper subset of \( X \), \( 0 \in S \), \( F:
S
\rightarrow R \)  be a functional on \( S \), \( F(0)=0
\),
\( p
\) be a sub-linear functional on \( X \), and
\begin{eqnarray*}
F(s_{1})+F(s_{2}) \leq p(s_{1}+s_{2}),s_{1},s_{2} \in S.
\end{eqnarray*}
Then \( F \) has an extension \( \hat{F}: X
\rightarrow R \) satisfying
\begin{eqnarray*}
\hat{F}(s)=F(s), s \in S,
\end{eqnarray*}
and
\begin{eqnarray*}
\hat{F}(x_{1})+\hat{F}(x_{2}) \leq p(x_{1}+x_{2}),x_{1},x_{2} \in X.
\end{eqnarray*}

\end{corollary}

\textsc{Proof.} \( S \) can be assumed to have at least two elements since
\( F \) can be extended to have the domain including the zero element and the
other element and to satisfy the required inequality otherwise. By Theorem 4,
there exists an extension
\(
\hat{F}
\) of
\( F
\) such that
\begin{eqnarray*}
\hat{F}(x_{1})+\hat{F}(x_{2}) \leq p(x_{1}+x_{2}),  x_{1}, x_{2} \in X, x_{1} \neq x_{2}.
\end{eqnarray*}
Further, because \( p \) is sub-linear,
\begin{eqnarray*}
\hat{F}(x)+\hat{F}(x) \leq p(x)+p(x)=p(x+x)
\end{eqnarray*}
for \( x \in X \).

\begin{flushright}
$\diamondsuit$
\end{flushright}

Let \( F^{+} \) and
\( F^{-} \) be the positive and negative parts of \( F: S \rightarrow R \) defined by
\( F^{+}(s)=\max[F(s),0] \) and \( F^{-}(s)=\max[-F(s),0] \) for \( s \in S \), respectively,
where \( S \) is a subset of the vector space \( X \). By the nonlinear Hahn-Banach theorem,
there exists a bounded extension to
the whole space for the functional having the bounded positive and negative parts
on the subset of the vector space.

\begin{corollary}

Let \( F \in B(S,R) \), \( 0 \in S \), \( F(0)=0 \), and for \( s_{1},s_{2} \in S \),
\begin{eqnarray*}
F^{+}(s_{1})+F^{+}(s_{2}) \leq M_{1}\left|\left|s_{1}+s_{2}\right|\right|_{X},
\end{eqnarray*}
and
\begin{eqnarray*}
F^{-}(s_{1})+F^{-}(s_{2}) \leq M_{2}\left|\left|s_{1}+s_{2}\right|\right|_{X},
\end{eqnarray*}
where \( S \) is a proper subset of a normed space \( X \), \( M_{1} \) and \( M_{2} \) are some positive constants, and \( F^{+} \) and
\( F^{-} \) are the positive and negative parts of \( F \), respectively.
Then there exists an extension \( \hat{F} \in B(X,R) \) of \( F \) satisfying
\begin{eqnarray*}
\hat{F}(s)=F(s), s \in S,
\end{eqnarray*}
and the inequality
\begin{eqnarray*}
\left|\hat{F}(x_{1})+\hat{F}(x_{2})\right| \leq
\left(M_{1}+M_{2}\right)\left|\left|x_{1}+x_{2}\right|\right|_{X}
\end{eqnarray*}
for \( x_{1},x_{2} \in X \).

\end{corollary}

\textsc{Proof.} By Corollary 2, there exist \( \hat{F}^{+} \) and
\( \hat{F}^{-} \) such that \( \hat{F}^{+}(s)=F^{+}(s) \),
\( \hat{F}^{-}(s)=F^{-}(s) \), \( s \in S \),
\begin{eqnarray*}
\hat{F}^{+}(x_{1})+\hat{F}^{+}(x_{2}) \leq M_{1}\left|\left|x_{1}+
x_{2}\right|\right|_{X},
\end{eqnarray*}
and
\begin{eqnarray*}
\hat{F}^{-}(x_{1})+\hat{F}^{-}(x_{2}) \leq M_{2}\left|\left|x_{1}+x_{2}
\right|\right|_{X}
\end{eqnarray*}
for \( x_{1},x_{2} \in X \). Let \( \hat{F}=\hat{F}^{+}-\hat{F}^{-} \). Then
\( \hat{F}(s)=F(s) \) and
\begin{eqnarray*}
&& \left|\hat{F}(x_{1})+\hat{F}(x_{2})\right| \\
&=&\left|\hat{F}^{+}(x_{1})+\hat{F}^{+}(x_{2})-
\hat{F}^{-}(x_{1})-\hat{F}^{-}(x_{2})\right| \\
&\leq& \left[\hat{F}^{+}(x_{1})+\hat{F}^{+}(x_{2})\right]+
\left[\hat{F}^{-}(x_{1})+\hat{F}^{-}(x_{2})\right]\\
 &\leq& \left(M_{1}+M_{2}\right)\left|\left|x_{1}+x_{2}\right|\right|_{X}.
 \end{eqnarray*}
Finally, \( \hat{F} \in B(X,R) \) because
\begin{eqnarray*}
\left|\hat{F}(x_{1})\right| \leq \left(M_{1}+M_{2}\right)\left|\left|x_{1}
\right|\right|_{X}.
 \end{eqnarray*}

\begin{flushright}
$\diamondsuit$
\end{flushright}

The complex version of the nonlinear Hahn-Banach theorem can be established
by applying Theorem 4, as stated by the following corollary. \\

\begin{corollary}

Let \( F: S \rightarrow C \) and \( F=F_{r}+iF_{c} \), where
\( S \) is a proper subset of a vector space \( X \) and both
\( F_{r} \) and \( F_{c} \) are real-valued functionals defined on \( S \). If for \( s_{1}, s_{2} \in S \)
and \( s_{1} \neq s_{2} \),
\begin{eqnarray*}
F_{r}(s_{1})+F_{r}(s_{2}) \leq p_{r}(s_{1}+s_{2}),
\end{eqnarray*}
and
\begin{eqnarray*}
F_{c}(s_{1})+F_{c}(s_{2}) \leq p_{c}(s_{1}+s_{2}),
\end{eqnarray*}
then \( F \) has an extension \( \hat{F}: X \rightarrow C \), \(
\hat{F}=\hat{F}_{r}+i\hat{F}_{c} \) satisfying
\begin{eqnarray*}
\hat{F}(s)=F(s),~s \in S,
\end{eqnarray*}
and for \( x_{1}, x_{2} \in X \) and \( x_{1} \neq x_{2} \),
\begin{eqnarray*}
\hat{F}_{r}(x_{1})+\hat{F}_{r}(x_{2}) \leq p_{r}(x_{1}+x_{2}),
\end{eqnarray*}
and
\begin{eqnarray*}
\hat{F}_{c}(x_{1})+\hat{F}_{c}(x_{2}) \leq p_{c}(x_{1}+x_{2}),
\end{eqnarray*}
where both \( \hat{F}_{r} \) and \( \hat{F}_{c} \) are real-valued functionals
defined on
\( X \) and
\( p_{r} \) and \( p_{c}
\) are sub-additive functionals on \( X \).

If \( 0 \in S \), \( F(0) =0 \),
\( p_{r} \) and \( p_{c} \) are sub-linear functionals on \( X \), and the
above inequalities for \( F_{r} \) and \( F_{c} \) hold for any \( s_{1},s_{2} \in S
\), then the above extension result holds and the above inequalities for
\( \hat{F}_{r} \) and \( \hat{F}_{c} \) also
hold for any \( x_{1},x_{2} \in X \).

\end{corollary}

The following theorem indicates that both the continuity and the inequality can
be preserved as extending a nonlinear functional on a closed subspace of a
separable Hilbert space to the whole space.

\begin{theorem}

Let \( Z \) with the orthonormal basis \( \{e_{j}\} \) be a proper closed subspace
of a separable Hilbert space
\( X \) with the orthonormal basis \( \{e_{j}\} \cup \{e^{\ast}_{i}\} \),
\( i=1,\ldots,j=1,\ldots \),  and
\( F:Z \rightarrow R \) be a continuous functional satisfying
\( F(0)=0 \) and
\begin{eqnarray*}
F(z_{1})+F(z_{2}) \leq p(z_{1}+z_{2})
\end{eqnarray*}
for orthogonal vectors \( z_{1} \) and \( z_{2} \) falling in the union of
one-dimensional spaces each spanned by \( e_{j}
\), where
\( p \) is a uniformly continuous sub-additive functional on
\( X \) with \( p(0)=0 \). Then there exists an extension \( \hat{F} \) of \( F \)
such that
\( \hat{F}: X \rightarrow R \) is continuous,
\begin{eqnarray*}
\hat{F}(z)=F(z), z \in Z,
\end{eqnarray*}
and
\begin{eqnarray*}
\hat{F}(x_{1})+\hat{F}(x_{2}) \leq p(x_{1}+x_{2})
\end{eqnarray*}
for orthogonal vectors \( x_{1} \) and \( x_{2} \) falling in the union of
one-dimensional spaces each spanned by \( e_{j} \) or \( e^{\ast}_{i} \).

\end{theorem}

\textsc{Proof.} Similar to the proof of Theorem 4, the existence of the maximal
element \( \hat{F} \) can be proved by Zorn's lemma and by defining the partial
ordering as the inclusion of the domains \( X_{m} \), \( m=1,\ldots
\), the space spanned by
\( \{e_{j}\} \) and \( \{e^{\ast}_{i}: i \leq m\} \), or \( X_{0} \) of the extensions, where
\( X_{0}=Z \). It remains to prove that
\( \hat{F} \) is defined on the whole space \( X \).

Suppose that \( Dom(\hat{F})=X_{m-1} \) is a proper closed subspace of \( X \). Let
\( E_{m-1} \subset X_{m-1} \) be the union of one-dimensional spaces each spanned by
\( \{e_{j}\} \) and  \( \{e^{\ast}_{i}: i \leq m-1 \} \), \( e^{\ast}_{0}=0 \), and
\( F_{1}:X_{m} \rightarrow R \) defined by
\begin{eqnarray*}
F_{1}(x+te^{\ast}_{m})=\hat{F}(x)-r(t),
\end{eqnarray*}
where \( x \in Dom(\hat{F}) \), \( t \in K \), and the real-valued function
\( r: K \rightarrow R \) is defined by
\begin{eqnarray*}
r(t)=\sup_{x \in E_{m-1}} \left[\hat{F}(x)-p(x+te^{\ast}_{m})\right].
\end{eqnarray*}
Then, for \( x \in E_{m-1} \),
\begin{eqnarray*}
\hat{F}(x)-p(x+te^{\ast}_{m}) \leq r(t)
\end{eqnarray*}
and hence
\begin{eqnarray*}
F_{1}(x)+F_{1}(te^{\ast}_{m})=\hat{F}(x)-r(t) \leq p(x+te^{\ast}_{m}).
\end{eqnarray*}

If \( F_{1} \) is continuous, then
\( F_{1} \) satisfying the required inequality is an extension of
the continuous functional \( \hat{F} \), i.e., a contradiction. It remains to prove
the continuity of \( F_{1} \). For
\( z_{n}=x_{n}+t_{n}e^{\ast}_{m}, x_{n} \in Dom(\hat{F}) \),
\( z_{n} \mathop{\longrightarrow}\limits_
{n \rightarrow\infty} z=x+te^{\ast}_{m}, x \in Dom(\hat{F}) \), implies that
\( x_{n} \mathop{\longrightarrow}\limits_
{n \rightarrow\infty} x \) and
\( t_{n} \mathop{\longrightarrow}\limits_
{n \rightarrow\infty} t \) owing to
\begin{eqnarray*}
\left|\left|z_{n}-z\right|\right|^{2}_{X}=
\left|\left|x_{n}-x\right|\right|^{2}_{X}+
\left|t_{n}-t\right|^{2},
\end{eqnarray*}
where the norm \( ||x||_{X}=(<x,x>_{X})^{1/2} \) and \( <\cdot,\cdot>_{X} \) is
the inner product on \( X \). Hence, if
\( r
\) is a continuous function of \( t \), then
\begin{eqnarray*}
\left|F_{1}(z_{n})-F_{1}(z)\right| \leq
\left|\hat{F}(x_{n})-\hat{F}(x)\right|+
\left|r(t_{n})-r(t)\right|
\end{eqnarray*}
and thus \( F_{1}(z_{n}) \mathop{\longrightarrow}\limits_ {n \rightarrow\infty}
F_{1}(z) \) as \( z_{n} \mathop{\longrightarrow}\limits_ {n \rightarrow\infty} z
\), i.e., \( F_{1} \) being continuous. For every \( \epsilon > 0 \), there exists an \( N \)
such that
\begin{eqnarray*}
&& r(t_{n}) \\ &=& \sup_{x \in E_{m-1}}\
\left[\hat{F}(x)-p(x+te^{\ast}_{m})+p(x+te^{\ast}_{m})-p(x+t_{n}e^{\ast}_{m})\right] \\
&\leq& \sup_{x \in E_{m-1}}\
\left[\hat{F}(x)-p(x+te^{\ast}_{m})\right]+ \sup_{x \in E_{m-1}}\
\left[p(x+te^{\ast}_{m})-p(x+t_{n}e^{\ast}_{m})\right] \\
&\leq& r(t)+ \epsilon
\end{eqnarray*}
for \( n \geq N \) by the uniform continuity of \( p \) and similarly
\( r(t) \leq r(t_{n})+\epsilon \),  i.e.,
\( |r(t)-r(t_{n})| \leq \epsilon \), and hence \( r \) is continuous.

\begin{flushright}
$\diamondsuit$
\end{flushright}

By Theorem 5, the continuity of the nonlinear functional and the boundedness on the basis can be
extended from the subspace to the whole space, as indicated by the following
corollary.

\begin{corollary}

Let \( Z \) with the orthonormal basis \( \{e_{j}\} \) be a proper closed subspace
of a separable Hilbert space
\( X \) with the orthonormal basis \( \{e_{j}\} \cup \{e^{\ast}_{i}\} \),
\( i=1,\ldots,j=1,\ldots \) and the norm induced by the inner product,
\( F:Z \rightarrow R \) be a continuous functional,
\( F(0)=0 \), and for orthogonal vectors \( z_{1} \) and \( z_{2} \) falling
in the union of one-dimensional spaces each spanned by \( e_{j}
\),
\begin{eqnarray*}
F^{+}(z_{1})+F^{+}(z_{2}) \leq M_{1}\left|\left|z_{1}+z_{2}\right|\right|_{X},
\end{eqnarray*}
and
\begin{eqnarray*}
F^{-}(z_{1})+F^{-}(z_{2}) \leq M_{2}\left|\left|z_{1}+z_{2}\right|\right|_{X},
\end{eqnarray*}
where \( M_{1} \) and \( M_{2} \) are some positive constants and \( F^{+} \) and
\( F^{-} \) are the positive and negative parts of \( F \), respectively. Then
there exists an extension \( \hat{F}: X \rightarrow R \) of \( F \) which is
continuous and satisfies
\begin{eqnarray*}
\hat{F}(z)=F(z), z \in Z,
\end{eqnarray*}
and the inequality
\begin{eqnarray*}
\left|\hat{F}(x_{1})+\hat{F}(x_{2})\right| \leq
\left(M_{1}+M_{2}\right)\left|\left|x_{1}+x_{2}\right|\right|_{X}
\end{eqnarray*}
for orthogonal vectors \( x_{1} \) and \( x_{2} \) falling in the union of
one-dimensional spaces each spanned by \( e_{j} \) or \( e^{\ast}_{i} \).

\end{corollary}

The complex version of Theorem 5 is stated by the following corollary. \\

\begin{corollary}

Let \( Z \) with the orthonormal basis \( \{e_{j}\} \) be a proper closed subspace
of a separable Hilbert space
\( X \) with the orthonormal basis \( \{e_{j}\} \cup \{e^{\ast}_{i}\} \),
\( i=1,\ldots,j=1,\ldots \), \( F: Z \rightarrow C \), \( F(0)=0 \), \( F=F_{r}+iF_{c} \), and
for orthogonal vectors \( z_{1} \) and \( z_{2} \) falling in the union of
one-dimensional spaces each spanned by \( e_{j} \),
\begin{eqnarray*}
F_{r}(z_{1})+F_{r}(z_{2}) \leq p_{r}(z_{1}+z_{2}),
\end{eqnarray*}
and
\begin{eqnarray*}
F_{c}(z_{1})+F_{c}(z_{2}) \leq p_{c}(z_{1}+z_{2}),
\end{eqnarray*}
where both
\( F_{r} \) and \( F_{c} \) are real-valued continuous functionals
defined on \( Z \) and \( p_{r} \) and \( p_{c} \) are uniformly continuous
sub-additive functionals on
\( X \) with \( p_{r}(0)=0 \) and \( p_{c}(0)=0 \). Then there exists an extension
\( \hat{F}: X \rightarrow C \)
such that
\begin{eqnarray*}
\hat{F}(z)=F(z), z \in Z,
\end{eqnarray*}
\(
\hat{F}=\hat{F}_{r}+i\hat{F}_{c}\)
and for orthogonal vectors \( x_{1} \) and \( x_{2} \) falling in the union of
one-dimensional spaces each spanned by \( e_{j} \) or
\( e^{\ast}_{i} \),
\begin{eqnarray*}
\hat{F}_{r}(x_{1})+\hat{F}_{r}(x_{2}) \leq p_{r}(x_{1}+x_{2}),
\end{eqnarray*}
and
\begin{eqnarray*}
\hat{F}_{c}(x_{1})+\hat{F}_{c}(x_{2}) \leq p_{c}(x_{1}+x_{2}),
\end{eqnarray*}
where both \( \hat{F}_{r} \) and \( \hat{F}_{c} \) are real-valued continuous
functionals defined on
\( X \).

\end{corollary}

The following theorem is a direct application of classic Hahn-Banach theorem
to the nonlinear case. The specific form of the nonlinear functional can be
preserved as extending from a subspace of a vector space to the whole space.
Let \( R^{+}=\{x: x \geq 0, x \in R\} \).

\begin{theorem}

Let \( F: Z \rightarrow R \) be a functional defined by
\begin{eqnarray*}
F(z)=f[|T(z)|]
\end{eqnarray*}
and
\begin{eqnarray*}
F(z) \leq f[p(z)], z \in Z,
\end{eqnarray*}
where
\( Z \) is a proper subspace of a vector space \( X \), \( T:Z \rightarrow K \)
is a linear functional, \( p \) is a semi-norm defined on \( X
\), and
\( f: R \rightarrow R \) is increasing on
\( R^{+} \). Then \( F \) has an extension
\( \hat{F}: X \rightarrow R \) with the form
\begin{eqnarray*}
\hat{F}(x)=f[|\hat{T}(x)|]
\end{eqnarray*}
and
\begin{eqnarray*}
\hat{F}(x) \leq f[p(x)],  x \in X,
\end{eqnarray*}
where
\( \hat{T}: X \rightarrow K  \) is a linear extension of \( T \), i.e., \( \hat{T}(z)=T(z) \) for
\( z \in Z \).

\end{theorem}

\textsc{Proof.} Because \( f \) is increasing on \( R^{+} \), then
\( F(z)=f[|T(z)|] \leq f[p(z)] \) implies that \( |T(z)| \leq p(z) \) for \( z \in Z \).
By the Hahn-Banach extension theorem (see \cite{rudin}, Theorem 3.3) there
exists a linear extension \( \hat{T} \) of \( T \) such that
\( \hat{T}(z)=T(z) \) for \( z \in Z \) and \( |\hat{T}(x)| \leq p(x) \) for \( x \in X \). Therefore,
let \( \hat{F}(x)=f[|\hat{T}(x)|] \). Then
\( \hat{F}(z)=F(z) \) for \( z \in Z \) and \(
\hat{F}(x)=f[|\hat{T}(x)|] \leq f[p(x)] \)  for \( x \in X \) owing to
\( f \) being increasing on \( R^{+} \).

\begin{flushright}
$\diamondsuit$
\end{flushright}

The above theorem can be applied to the nonlinear functionals associated with the
powers of the linear functional.

\begin{corollary}

Let \( F:Z \rightarrow R \) be a functional defined by
\begin{eqnarray*}
F(z)=|T(z)|^{k}
\end{eqnarray*}
and
\begin{eqnarray*}
F(z) \leq [p(z)]^{k}, z \in Z,
\end{eqnarray*}
where \( k \) is a positive number, \( Z \) is a proper subspace of a vector space
\( X
\),
\( T:Z \rightarrow K \) is a linear functional, and
\( p \) is a semi-norm defined on \( X \). Then \( F \) has an extension
\( \hat{F}: X \rightarrow R \) with the form
\begin{eqnarray*}
\hat{F}(x)=|\hat{T}(x)|^{k}
\end{eqnarray*}
and
\begin{eqnarray*}
\hat{F}(x) \leq [p(x)]^{k}, x \in X,
\end{eqnarray*}
where
\( \hat{T}: X \rightarrow K \) is a linear extension of \( T \).

\end{corollary}

\textsc{Proof.} Let  \( f(x)=x^{k} \) and hence the results hold by Theorem 6.

\begin{flushright}
$\diamondsuit$
\end{flushright}

The above nonlinear functionals associated with the bounded linear functionals
can be bounded as certain norms are employed. As
the linear functional \( T \) is bounded, the
following corollary indicates that the functionals in Corollary 7 fall in
\( [B(Z,R),p_{k}] \) (see Remark 1) and the associated extended functionals are in \(
[B(X,R),p_{k}] \).

\begin{corollary}

Let \( F:Z \rightarrow R \) be a functional defined by
\begin{eqnarray*}
F(z)=|T(z)|^{k}, z \in Z,
\end{eqnarray*}
where \( k \) is a positive number, \( Z \) is a proper subspace of a normed space
\( X \), and \( T:Z
\rightarrow K
\) is a bounded linear functional. Then \( F
\in [B(Z,R),p_{k}] \) and there exists an extension \( \hat{F}:X \rightarrow R \) of
\( F \) such that \( \hat{F} \in [B(X,R),p_{k}]  \),
\begin{eqnarray*}
||\hat{F}||_{[B(X,R),p_{k}]}=\left|\left|F\right|\right|_{[B(Z,R),p_{k}]},
\end{eqnarray*}
and
\begin{eqnarray*}
\hat{F}(x)=|\hat{T}(x)|^{k}, x \in X,
\end{eqnarray*}
where the bounded linear functional
\( \hat{T}:X \rightarrow K \) is an extension of \( T \).

\end{corollary}

\textsc{Proof.}
\begin{eqnarray*}
\left|\left|F\right|\right|_{\left[B(Z,R),p_{k}\right]}=
\max\left[
\sup_{z \neq 0, z \in Z}\frac{\left|F(z)\right|}{\left|\left|z\right|\right|^{k}_{X}},
\left|F(0)\right|\right]= \left|\left|T\right|\right|_{B(Z,K)}^{k}
\end{eqnarray*}
is finite, i.e., \( F \in [B(Z,R),p_{k}] \). Furthermore, by the Hahn-Banach
theorem, there exists a linear extension
\( \hat{T}:X \rightarrow K \) of \( T \) such that \( ||T||_{B(Z,K)}=||\hat{T}||_{B(X,K)} \).
Let \( \hat{F}(x)=|\hat{T}(x)|^{k}, x \in X \), then \( \hat{F}(z)=F(z), z \in Z, \) and
\begin{eqnarray*}
\left|\left|\hat{F}\right|\right|_{\left[B(X,R),p_{k}\right]}=
\left(
\sup_{x \neq 0, x \in X}\frac{\left|\hat{T}(x)\right|}{\left|\left|x
\right|\right|_{X}}\right)^{k} =\left|\left|\hat{T}\right|\right|^{k}_{B(X,K)}
=\left|\left|F\right|\right|_{\left[B(Z,R),p_{k}\right]},
\end{eqnarray*}
i.e., \( \hat{F} \in [B(X,R),p_{k}] \).

\begin{flushright}
$\diamondsuit$
\end{flushright}

The following corollary is the complex version of Theorem 6. \\

\begin{corollary}

Let \( F: Z \rightarrow C \) and \( F=F_{r}+iF_{c} \), where both \( F_{r} \) and
\( F_{c} \) are real-valued functionals on \( Z  \) defined by
\begin{eqnarray*}
F_{r}(z)=f_{r}[|T_{r}(z)|]
\end{eqnarray*}
and
\begin{eqnarray*}
F_{c}(z)=f_{c}[|T_{c}(z)|]
\end{eqnarray*}
and satisfying
\begin{eqnarray*}
F_{r}(z) \leq f_{r}[p_{r}(z)]
\end{eqnarray*}
and
\begin{eqnarray*}
F_{c}(z) \leq f_{c}[p_{c}(z)]
\end{eqnarray*}
for \( z \in Z \), and where
\( Z \) is a proper subspace of a vector space \( X \),
\( T_{r}:Z \rightarrow K \) and \( T_{c}:Z \rightarrow K \) are linear functionals
defined on \( Z \), both \( p_{r} \) and \( p_{c} \) are semi-norms defined on
\( X \), and
\( f_{r}: R \rightarrow R \) and \( f_{c}: R
\rightarrow R \) are  increasing on
\( R^{+} \). Then \( F \) has an extension
\( \hat{F}: X \rightarrow C \) having the form
\( \hat{F}=\hat{F}_{r}+i\hat{F}_{c} \),
\begin{eqnarray*}
\hat{F}_{r}(x)=f_{r}[|\hat{T}_{r}(x)|]
\end{eqnarray*}
and
\begin{eqnarray*}
\hat{F}_{c}(x)=f_{c}[|\hat{T}_{c}(x)|]
\end{eqnarray*}
and satisfying
\begin{eqnarray*}
\hat{F}_{r}(x) \leq f_{r}[p_{r}(x)]
\end{eqnarray*}
and
\begin{eqnarray*}
\hat{F}_{c}(x) \leq f_{c}[p_{c}(x)]
\end{eqnarray*}
for \( x \in X \), where both \( \hat{F}_{r} \) and
\( \hat{F}_{c} \) are real-valued functionals defined on \( X \) and
\( \hat{T}_{r}:X \rightarrow K \) and \( \hat{T}_{c}:X \rightarrow K \)  are linear extensions of
\( T_{r} \) and \( T_{c} \), respectively.

\end{corollary}

\section{Nonlinear mappings on metrizable spaces}

The operators of interest in previous sections are defined on the vector spaces or the normed
spaces. In this section, the mappings on a
metrizable space can have normed values by defining a norm. Thus, some
linear functional such as the Dirac's delta function considered as the linear
functional falls in certain normed spaces. The basic facts about the mappings of
interest are given in next subsection, while several examples of the
mappings on the metrizable spaces are presented in the second subsection.

\subsection{Bounded mappings on metrizable spaces}

In this subsection, the main results that
the set of mappings between certain metrizable spaces
being a translation invariant metrizable vector space and some mappings from a metrizable
space into a normed space falling in a normed space are proved in Theorem 10.
Thus some linear functionals corresponding to commonly used
distributions fall in some normed spaces and are given in next subsection.

In this subsection, let \( X \) and \( Y \) be the metrizable vector spaces over the
field
\( K \). Let
\( V(X,Y) \) be the vector space of all mappings from \( X  \) into \( Y \) with
the algebraic operations of \( F_{1}, F_{2} \in V(X,Y) \) being the mappings from
\( X \) into
\( Y \) defined by \( (F_{1}+F_{2})(x)=F_{1}(x)+F_{2}(x) \) and
\( (\alpha F_{1})(x)=\alpha F_{1}(x) \) for \( x \in X \) and \( \alpha \in K \).
Note that the zero element in
\( V(X,Y) \) is the mapping of which image equal to the zero element in \( Y \).
Define a nonnegative extended real-valued function \( d \) on \( V(X,Y) \) by
\begin{eqnarray*}
d(F_{1},F_{2})=\max\left\{\sup_{x \neq 0,x \in X}
\frac{d_{Y}\left[F_{1}(x),F_{2}(x)\right]}{d_{X}(x,0)},
d_{Y}\left[F_{1}(0),F_{2}(0)\right]\right\},
\end{eqnarray*}
where \( d_{X} \) and \( d_{Y} \) are the metrics on
\( X \) and \( Y \), respectively. Let \( B_{d}(X,Y) \), containing the zero element of
\(  V(X,Y) \) and the
subset of \( V(X,Y) \), have the property that \( d(F_{1},F_{2}) < \infty \)
for any \( F_{1}, F_{2} \in B_{d}(X,Y) \). If \( F \) is a mapping from
\( X \) into \( Y \) and \( d(F,0) \) is finite, then \( F \in B_{d}(X,Y) \) because
for any \( G \in B_{d}(X,Y) \),
\begin{eqnarray*}
d(F,G) \leq d(F,0)+d(G,0)
\end{eqnarray*}
is finite.

The following theorem indicates that \( [B_{d}(X,Y),d] \) is a metric space. The routine proof is not presented.

\begin{theorem}

\( d \) is a metric on \( B_{d}(X,Y) \) and \( [B_{d}(X,Y),d] \) is a metric space.

\end{theorem}

It is well known that the space of all bounded linear operators from
a normed space X to a Banach space Y is complete. The following
corollary can be considered as the generalization of the completeness
result for the bounded linear operators to the possibly nonlinear mappings
on the metrizable spaces.

\begin{corollary}

If \( Y \) is complete, then \( [B_{d}(X,Y),d] \) is a complete metric space.

\end{corollary}

\textsc{Proof.} Let \( \{F_{n}\} \) be a Cauchy sequence in \( B_{d}(X,Y) \).
Then for any positive \( \epsilon \), there exists an \( N \) such that for \( m, n
> N \), \( d(F_{n},F_{m}) < \epsilon \). Then as \( x \neq 0 \),
\begin{eqnarray*}
d_{Y}\left[F_{n}(x),F_{m}(x)\right] \leq d(F_{n},F_{m})d_{X}(x,0) <
\epsilon d_{X}(x,0)
\end{eqnarray*}
and
\begin{eqnarray*}
d_{Y}\left[F_{n}(0),F_{m}(0)\right] \leq d(F_{n},F_{m}) < \epsilon.
\end{eqnarray*}
Thus, \( \{F_{n}(x)\} \) is Cauchy in \( Y \) for \( x \in X \) and
\( F_{n}(x) \mathop{\longrightarrow}\limits_ {n
\rightarrow\infty} y, y \in Y \) owing to the completeness of \( Y \). Define a mapping
\( F: X \rightarrow Y \) by \( F(x)=y \). For
\( x \neq 0 \),
\begin{eqnarray*}
d_{Y}\left[F_{n}(x),F(x)\right] &=& d_{Y}\left[F_{n}(x),\lim_{m \rightarrow
\infty} F_{m}(x)\right]=\lim_{m \rightarrow
\infty}d_{Y}\left[F_{n}(x),F_{m}(x)\right] \\ &\leq &\epsilon d_{X}(x,0).
\end{eqnarray*}
In addition,
\begin{eqnarray*}
d_{Y}\left[F_{n}(0),F(0)\right] = \lim_{m \rightarrow \infty}
d_{Y}\left[F_{n}(0),F_{m}(0)\right] \leq \epsilon.
\end{eqnarray*}
Thus, \( F \in [B_{d}(X,Y),d] \) owing to \( d(F,F_{n}) \) being finite and \( \{F_{
n}\}
\) converges to
\( F
\) because of
\( d(F_{n},F)
\leq
\epsilon.
\)

\begin{flushright}
$\diamondsuit$
\end{flushright}

The following theorem gives the characterization of the mappings falling in
\( B_{d}(X,Y) \).

\begin{theorem}

If \( F \in B_{d}(X,Y) \), then for any
\( \epsilon_{1} > 0 \) there exists a \( \epsilon_{2} > 0 \) such that
\( F[B_{X}(\epsilon_{1})] \subset B_{Y}(\epsilon_{2}) \), where \( B_{X}(\epsilon_{1})=\{x: d_{X}(x,0) < \epsilon_{1},x \in X\}
\) and \( B_{Y}(\epsilon_{2})=\{y: d_{Y}(y,0) < \epsilon_{2},y \in Y\}
\). On the other hand, if the limit
\begin{eqnarray*}
\lim_{x \rightarrow 0} \frac{d_{Y}\left[F(x),0\right]}{d_{X}(x,0)}
\end{eqnarray*}
exists and is finite, there exists a bounded ball \( B_{X}(\epsilon) \) such that \( F
\) maps the complement of
\( B_{X}(\epsilon) \) into a bounded subset of some bounded ball \( B_{Y}(\epsilon^{\ast})
\), and  for any
\( \epsilon_{1} > 0 \) there exists a \( \epsilon_{2} > 0 \) such that
\( F[B_{X}(\epsilon_{1})] \subset B_{Y}(\epsilon_{2}) \), then
\( F \in B_{d}(X,Y) \), where \( \epsilon, \epsilon^{\ast} > 0 \).

\end{theorem}

\textsc{Proof.} If \( F \in B_{d}(X,Y) \), then \( d(F,0)  \) is finite. Hence for any \( x \in X \),
\begin{eqnarray*}
d_{Y}\left[F(x),0\right] \leq d(F,0)
\max\left[d_{X}(x,0),1\right]
\end{eqnarray*}
and thus
\begin{eqnarray*}
d_{Y}\left[F(x^{\ast}),0\right] \leq d(F,0) \max\left(\epsilon_{1},1\right)
\end{eqnarray*}
for any \( x^{\ast} \in B_{X}(\epsilon_{1}) \).

On the other hand, if the given conditions hold, there exist positive numbers
\( \delta < \epsilon  \), \( m \), and \( M \) such that for
\( d_{X}(x_{1},0) < \delta,  x_{1} \neq 0 \), and \( d_{X}(x_{2},0) > \epsilon \),
\begin{eqnarray*}
\frac{d_{Y}\left[F(x_{1}),0\right]}{d_{X}(x_{1},0)} < m
\end{eqnarray*}
and
\begin{eqnarray*}
\frac{d_{Y}\left[F(x_{2}),0\right]}{d_{X}(x_{2},0)} < M.
\end{eqnarray*}
Therefore,
\begin{eqnarray*}
&& d(F,0) \\ &=& \max\left\{
\sup_{x \neq 0, x \in X}\frac{d_{Y}\left[F(x),0\right]}{d_{X}(x,0)},
d_{Y}\left[F(0),0\right]\right\} \\ &=& \max\left\{
\sup_{d_{X}(x,0) < \delta, x \neq 0, x \in X}\frac{d_{Y}\left[F(x),0\right]}
{d_{X}(x,0)},
\sup_{\delta \leq d_{X}(x,0) \leq \epsilon, x \in X}
\frac{d_{Y}\left[F(x),0\right]}{d_{X}(x,0)}, \right. \\
&& \left. \sup_{d_{X}(x,0) > \epsilon, x \in X}\frac{d_{Y}\left[F(x),0\right]}
{d_{X}(x,0)},d_{Y}\left[F(0),0\right]\right\}\\ &\leq&
\max\left\{m,\frac{\overline{\epsilon}}{\delta},M, d_{Y}\left[F(0),0\right]\right\},
\end{eqnarray*}
and thus \( F \in B_{d}(X,Y) \), where \( \overline{\epsilon} \) is some positive number.

\begin{flushright}
$\diamondsuit$
\end{flushright}

It is natural to ask when \( B_{d}(X,Y) \) can be a vector space. It turns out that
the following property plays a crucial role.

\begin{definition}

A translation invariant metric \( d_{X} \) is scale bounded on a subset \( S \) of a
metric vector space
\( [X,d_{X}]
\) if and only if
\begin{eqnarray*}
d_{X}(\alpha s,0) \leq C_{\alpha} d_{X}(s,0), s \in S,
\end{eqnarray*}
where \( \alpha \in K \) and \( C_{\alpha} \) is a positive number depending on \(
\alpha \).

\end{definition}

Note that the metric induced by the normed function is scale bounded.
A linear operator from a normed space into another normed space is continuous if and only
if it is bounded. In addition, if the linear operator is continuous at one point, then it is a bounded
operator. The continuity of a linear mapping implying the boundedness of the linear mapping
relies on the scale boundedness of the metric, as
indicated by the following theorem.

\begin{theorem}

Let \( F: X \rightarrow Y \) be a linear mapping and the metric \( d_{Y} \) is translation invariant, where
\( X \) and \( Y \) are metric vector spaces.  Then:  \\
 (a) If  \( F \in B_{d}(X,Y) \), then \( F \) is continuous. On the other hand, if
 \( F \) is continuous, the translation invariant metric \( d_{X} \) is scale bounded with \( C_{\alpha}=M(\alpha)|\alpha| \),
 \( Y \) is a normed space, and  the metric on \( Y \)  is the one induced by the norm, then \( F \in B_{d}(X,Y) \), where \( C_{\alpha} \) is considered
 as a positive function defined on \( K \) and \( M \) is a positive bounded function defined on \( K \). \\
 (b) If \( F \) is continuous at a single point and \( d_{X} \) is translation invariant, then \( F \) is continuous.

\end{theorem}

 \textsc{Proof.} (a): If \( F \in B_{d}(X,Y) \), then
\begin{eqnarray*}
d_{Y}\left[F(x_{n}),F(x)\right]=d_{Y}\left[F(x_{n}-x),0\right] \leq d(F,0)d_{X}(x_{n}-x,0)
\end{eqnarray*}
and thus  \( F(x_{n}) \mathop{\longrightarrow}\limits_ {n\rightarrow \infty}
F(x) \) as \( x_{n} \mathop{\longrightarrow}\limits_ {n \rightarrow \infty}
x \) in \( X \), i.e., \( F \)  being continuous. On the other hand, if \( F \) is continuous and
\( Y \) is a normed space with the metric induced by the norm, then for every \( \epsilon > 0 \)
there exists a \(\delta > 0 \) such that
 \( ||F(x)||_{Y} < \epsilon \) for \( x \in X \) satisfying \( d_{X}(x,0)< \delta \). Further,
 \( d_{X}[x/d_{X}(x,0), 0] \) is bounded for \( x \neq 0 \)  by the condition
 imposed on \( d_{X} \) and hence there exists an
 \( \alpha_{0} \neq 0  \) such that \( d[(\alpha_{0} x)/d_{X}(x,0), 0] < \delta \). Thus,
 for \( x \neq 0 \), \( ||F[(\alpha_{0} x)/d_{X}(x,0)]||_{Y}= |\alpha_{0}|||F[x/d_{X}(x,0)]||_{Y} < \epsilon \)
 and hence
 \(  ||F[x/d_{X}(x,0)]||_{Y} < \epsilon/ |\alpha_{0}| \). Finally,  \( F \in B_{d}(X,Y) \) because
 \begin{eqnarray*}
 d(F,0)=\sup_{x \neq 0,x \in X}
\frac{\left|\left|d_{X}(x,0)F\left[\frac{x}{d_{X}(x,0)}\right]\right|\right|_{Y}}{d_{X}(x,0)}
<  \epsilon/ |\alpha_{0}|.
\end{eqnarray*}
(b): Assume that \( F \) is continuous at \( x_{0} \). Thus, for every \( \epsilon > 0 \),
there exists a \( \delta > 0 \) such that
\( d_{Y}[F(x^{\ast}),F(x_{0})] < \epsilon \) for any \( x^{\ast} \in X \) satisfying
\( d_{X}(x^{\ast},x_{0}) < \delta \). Then, by the translation invariance property of the metrics,
for any \( x \in X \) and any \( x^{\ast\ast} \in X \) satisfying
\( d_{X}(x^{\ast\ast},x)=d_{X}(x_{0}+x^{\ast\ast}-x,x_{0}) < \delta \), \( d_{Y}[F(x^{\ast\ast}),F(x)]=
d_{Y}[F(x_{0}+x^{\ast\ast}-x),F(x_{0})]  < \epsilon \), i.e., \( F \) being continuous.

\begin{flushright}
$\diamondsuit$
\end{flushright}

The scale boundedness of the translation invariant metric \( d_{Y} \) on \( Y \) is the key
for \( B_{d}(X,Y) \) being a vector space.  Furthermore,
if \( Y \) is a normed space, it turns out that \( B_{d}(X,Y) \) can be a normed
space, as given in the following
theorem.

\begin{theorem}

If the translation invariant metric \( d_{Y} \) is scale bounded on \( Y \), then
\( [B_{d}(X,Y),d] \)  is a translation invariant metric vector space.
Furthermore, if \( Y \) is a normed space and \( d_{Y} \) is the metric induced by the norm on
\( Y \), then
\( B_{d}(X,Y)
\) is a normed space with the norm
\( ||F||_{B_{d}(X,Y)}=d(F,0) \).

\end{theorem}

\textsc{Proof.} To prove that \( d \) is translation invariant, for \( F_{1},F_{2},F_{3} \in B_{d}(X,Y)
\),
\begin{eqnarray*}
&& d(F_{1}+F_{3},F_{2}+F_{3}) \\ &=& \max\left\{\sup_{x \neq 0, x \in X}
\frac{d_{Y}
\left[F_{1}(x)+F_{3}(x),F_{2}(x)+F_{3}(x)\right]}{d_{X}(x,0)}, \right. \\
&& \left. d_{Y}
\left[F_{1}(0)+F_{3}(0),F_{2}(0)+F_{3}(0)\right]\right\} \\
&=& \max\left\{\sup_{x \neq 0, x \in X} \frac{d_{Y}
\left[F_{1}(x),F_{2}(x)\right]}{d_{X}(x,0)},d_{Y}
\left[F_{1}(0),F_{2}(0)\right]\right\} \\
&=&d(F_{1},F_{2})
\end{eqnarray*}
owing to \( d_{Y} \) being translation invariant. To prove that \( B_{d}(X,Y) \)
is a vector space, for \( F_{1},F_{2} \in B_{d}(X,Y) \) and \(
\alpha \in K \),
\begin{eqnarray*}
&& d(\alpha F_{1},0)\\ &=& \max\left\{\sup_{x \neq 0, x \in X} \frac{d_{Y}
\left[\alpha F_{1}(x),0)\right]}{d_{X}(x,0)},d_{Y}
\left[\alpha F_{1}(0),0\right]\right\} \\
&\leq& C_{\alpha} d(F_{1},0)
\end{eqnarray*}
is finite, where  \( C_{\alpha} \) is a positive number depending on \(
\alpha \). Then by the translation invariance
of \( d \),
\begin{eqnarray*}
d(\alpha F_{1}+F_{2},0) \leq d(\alpha F_{1},0)+d(F_{2},0)
\end{eqnarray*}
gives that \( d(\alpha F_{1}+F_{2},0) \) is finite, i.e.,
\( \alpha F_{1}+F_{2} \in B_{d}(X,Y) \) and \( B_{d}(X,Y) \) being a
vector space.

Next is to prove that \( B_{d}(X,Y)
\) is a normed space. Since \( Y \) is a normed space, \( d_{Y} \) induced by the norm
is translation invariant and scale bounded. Thus, \( B_{d}(X,Y) \) is a vector space. It remains to prove
that \( d(F,0) \) is a norm. For \( F_{1},F_{2} \in B_{d}(X,Y) \) and \(
\alpha \in K \), by the properties of the translation invariant metric
\( d
\),
\( ||F_{1}||_{B_{d}(X,Y)}=d(F_{1},0) \geq 0 \),
\( d(F_{1},0)=||F_{1}||_{B_{d}(X,Y)}=0 \) if and only if \( F_{1}=0 \),
\begin{eqnarray*}
&& \left|\left|F_{1}+F_{2}\right|\right|_{B_{d}(X,Y)} \\ &\leq&
d(F_{1}+F_{2},F_{2})+d(F_{2},0) \\ &=& d(F_{1},0)+d(F_{2},0) \\ &=&
\left|\left|F_{1}\right|\right|_{B_{d}(X,Y)}+\left|\left|F_{2}\right|\right|_{B_{d}(X,Y)},
\end{eqnarray*}
and finally
\begin{eqnarray*}
&& \left|\left|\alpha F_{1} \right|\right|_{B_{d}(X,Y)} \\ &=& \max\left\{
\sup_{x \neq 0, x \in X}\frac{\left|\left|\alpha F_{1}(x) \right|\right|_{Y}}{d_{X}(x,0)},
\left|\left|\alpha F_{1}(0)\right|\right|_{Y}\right\} \\
&=& \max\left\{
\sup_{x \neq 0, x \in X}\frac{\left|\alpha\right|\left|\left|F_{1}(x) \right|\right|_{Y}}{d_{X}(x,0)},
\left|\alpha\right|\left|\left|F_{1}(0)\right|\right|_{Y}\right\} \\
&=& \left|\alpha\right|\left|\left|F_{1} \right|\right|_{B_{d}(X,Y)}.
\end{eqnarray*}

\begin{flushright}
$\diamondsuit$
\end{flushright}

\begin{remark}

The space \( B(X,Y)
\) is one example of the space \( B_{d}(X,Y) \).
Let \( X \) and \( Y \) be the normed spaces. Then the space \( B(X,Y)
\) is the space \( B_{d}(X,Y) \) as
\( d_{x} \) and \( d_{y} \) are the metrics induces by the norms on
\( X \) and \( Y \), respectively, i.e., \( d \) being the metric induced by
the norm on \( B(X,Y) \).

\end{remark}

By Corollary 10 and Theorem 10, the following result holds.

\begin{corollary}

If \( Y \) is a Banach space and \( d_{Y} \) is the metric induced by the norm on
\( Y \), then \( B_{d}(X,Y) \) is a Banach space with the norm
\( ||F||_{B_{d}(X,Y)}=d(F,0) \).

\end{corollary}

The normed space of all bounded linear operators from a normed space \( X
\) into \( X \) is a normed algebra with the multiplication being the composition of the
operators. The normed space \( B(X) \) is not a normed algebra with the
multiplication being the composition of the operators. However,
\( B(X) \) can be a normed algebra or a Banach algebra depending on
\( X \) being a normed algebra or a Banach algebra as the multiplication
of the operators is defined properly. Similarly, the normed space \( B(X,Y) \)
can be a normed algebra or a Banach algebra depending on
\( Y \) being a normed algebra or a Banach algebra. The following corollary indicates that
\( B_{d}(X,Y) \) can be a normed algebra
 or a Banach algebra by defining an operation of
multiplication for two mappings and hence the normed space \( B(X,Y) \)
can be a normed algebra or a Banach algebra (see Remark 2).

\begin{corollary}

Let \( Y \) be a normed (Banach) algebra and \( d_{Y} \) is the metric induced by
the norm on
\( Y \). Define the multiplication of two mappings \( F_{1}: X \rightarrow Y \) and \( F_{2}: X \rightarrow Y \) by
\begin{eqnarray*}
\left(F_{1} \ast F_{2}\right)(x)=\frac{F_{1}(x)F_{2}(x)}{d_{X}(x,0)},x \neq 0, x \in X,
\end{eqnarray*}
and
\begin{eqnarray*}
\left(F_{1} \ast F_{2}\right)(0)=F_{1}(0)F_{2}(0).
\end{eqnarray*}
Then \( B_{d}(X,Y) \) is a normed (Banach) algebra with the norm
\( ||F||_{B_{d}(X,Y)}=d(F,0) \). If
\( Y \) is unital, then \( B_{d}(X,Y) \) is unital.

\end{corollary}

\textsc{Proof.} By Theorem 10  or Corollary 11, \( B_{d}(X,Y) \) is a normed
space or a Banach space depending on \( Y \) being a normed space or a Banach
space. Next is to prove that \( B_{d}(X,Y) \) is
a normed algebra. Let \( F_{1},F_{2}, F_{3} \in B_{d}(X,Y) \). First, as \( x \neq 0 \),
\begin{eqnarray*}
\left[(F_{1}*F_{2})*F_{3}\right](x)&=&\frac{(F_{1}*F_{2})(x)F_{3}(x)}
                                {d_{X}(x,0)} \\
                                 &=&\frac{F_{1}(x)F_{2}(x)F_{3}(x)}{\left[d_{X}(x,0)\right]^{2}} \\
                                 &=&\frac{F_{1}(x)(F_{2}*F_{3})(x)}{d_{X}(x,0)} \\
                                 &=&\left[F_{1}*(F_{2}*F_{3})\right](x).
\end{eqnarray*}
As \( x=0 \),
\begin{eqnarray*}
\left[(F_{1}*F_{2})*F_{3}\right](0)&=&(F_{1}*F_{2})(0)F_{3}(0)\\
                                 &=&F_{1}(0)F_{2}(0)F_{3}(0) \\
                                 &=&F_{1}(0)(F_{2}*F_{3})(0) \\
                                 &=&\left[F_{1}*(F_{2}*F_{3})\right](0).
\end{eqnarray*}
Secondly, as \( x \neq 0 \),
\begin{eqnarray*}
\left[F_{1}*(F_{2}+F_{3})\right](x)&=&\frac{F_{1}(x)(F_{2}+F_{3})(x)}
                                 {d_{X}(x,0)} \\
                                 &=&\frac{F_{1}(x)F_{2}(x)+F_{1}(x)F_{3}(x)}
                                 {d_{X}(x,0)} \\
                                 &=&(F_{1}*F_{2})(x)+(F_{1}*F_{3})(x).
\end{eqnarray*}
As \( x=0 \),
\begin{eqnarray*}
\left[F_{1}*(F_{2}+F_{3})\right](0)&=&F_{1}(0)(F_{2}+F_{3})(0)\\
                                 &=&F_{1}(0)F_{2}(0)+F_{1}(0)F_{3}(0)\\
                                 &=&(F_{1}*F_{2})(0)+(F_{1}*F_{3})(0).
\end{eqnarray*}
\( [(F_{1}+F_{2})*F_{3}](x)=(F_{1}*F_{3})(x)+(F_{2}*F_{3})(x) \) and
\( [\alpha(F_{1}*F_{2})](x)=
[(\alpha F_{1})*F_{2}](x)=[F_{1}*(\alpha F_{2})](x) \) for \( x \in X \) and
\( \alpha \in K \), can be
proved analogously. Further,
\begin{eqnarray*}
&& \|F_{1}*F_{2}\|_{B_{d}(X,Y)} \\ &=&\max\left(\sup_{x \neq 0, x \in
X}\frac{\|F_{1}(x)F_{2}(x)\|_{Y}}{\left[d_{X}(x,0)\right]^{2}},
\|F_{1}(0)F_{2}(0)\|_{Y}\right) \\
&\leq&\max\left(\sup_{x \neq 0, x \in X}\frac{\|F_{1}(x)\|_{Y}}{d_{X}(x,0)}
\sup_{x \neq 0, x \in X}\frac{\|F_{2}(x)\|_{Y}}{d_{X}(x,0)},
\|F_{1}(0)\|_{Y}\|F_{2}(0)\|_{Y}\right) \\
&\leq& \max\left(\sup_{x \neq 0, x \in X}\frac{\|F_{1}(x)\|_{Y}}
{d_{X}(x,0)},\|F_{1}(0)\|_{Y}\right) \\ && \max\left(\sup_{x \neq 0, x \in
X}\frac{\|F_{2}(x)\|_{Y}}{d_{X}(x,0)},
\|F_{2}(0)\|_{Y}\right) \\
&=& \|F_{1}\|_{B_{d}(X,Y)}\|F_{2}\|_{B_{d}(X,Y)}.
\end{eqnarray*}

As \( 1 \) is the unit element in \( Y \), the unit element \( e \) in \( B_{d}(X,Y) \) is
given by
\( e(x)=d_{X}(x,0) 1 \) for \( x \neq 0 \) and
\( e(0)=1 \). Then
\begin{eqnarray*}
(F*e)(x)&=&\frac{F(x)e(x)}{d_{X}(x,0)}=F(x)1=F(x) \\
       &=&1F(x)=\frac{e(x)F(x)}{d_{X}(x,0)}=(e*F)(x)
\end{eqnarray*}
for
\( x \neq 0 \),
\begin{eqnarray*}
(F*e)(0)&=&F(0)e(0)=F(0)1=F(0) \\
       &=&1F(0)=e(0)F(0)=(e*F)(0),
\end{eqnarray*}
and
\begin{eqnarray*}
\|e\|_{B_{d}(X,Y)}&=&\max\left(\sup_{x \neq 0, x \in X}\frac{\|e(x)\|_{Y}}{d_{X}(x,0)},
            \|e(0)\|_{Y}\right) \\
             &=& \|1\|_{Y} \\
	     &=& 1,
\end{eqnarray*}
where the last \( 1 \) is the unit element in the scalar field.

\begin{flushright}
$\diamondsuit$
\end{flushright}

\subsection{Examples}

The following examples are the applications of the results in previous
subsection. The first example is concerned with the linear functionals
corresponding to the distributions. By defining a translation invariant metric,
commonly used linear functionals such as the distributions corresponding to a
Lebesgue integrable function and the Dirac's delta function fall in some
Banach space by Corollary 11. The second example is concerned with the possibly
nonlinear operators defined on the Banach space which is one of the Banach
spaces given in the first example. The operators of interest are
associated with the position operator and the momentum operator in quantum mechanics.
Finally, the Fourier transform and the Fourier-Plancherel transform defined as the
bounded linear operators and the bounded linear mapping, respectively, are shown in the third
example.

\begin{example}

Let the subspace \( D(\Omega) \) of \( C^{\infty}(\Omega) \) be the space of all
complex valued functions defined on the nonempty open subset \( \Omega \)
of
\( R^{n}
\) with compact supports, where
\( C^{\infty}(\Omega) \) is the space of all complex valued infinitely
differentiable functions defined on \( \Omega \). Define a translation invariant
metric on \( C^{\infty}(\Omega) \) by
\begin{eqnarray*}
d_{C^{\infty}(\Omega) }(f,g)=
\max\left[\sum_{i=1}^{\infty}\frac{a^{-i}p_{i}(f-g)}{b+p_{i}(f-g)},
p_{N}(f-g)\right]
\end{eqnarray*}
for \( f, g \in C^{\infty}(\Omega) \), where \( N \) is a positive integer,
\( a
> 1
\),
\( b
> 0\),
\begin{eqnarray*}
p_{i}(f)=\max\left\{\left|\left(\frac{\partial^{|\alpha|} f} {\partial
x_{j_{1}}^{\alpha_{1}}
\cdots  \partial x_{j_{r}}^{\alpha_{r}}}\right)(x)\right|:
x=(x_{1},\ldots,x_{n})^{t} \in K_{i},
\left|\alpha\right| \leq i\right\},
\end{eqnarray*}
and where \( \alpha=(\alpha_{1},\ldots,\alpha_{r}) \),
\( |\alpha|=\sum_{k=1}^{r}\alpha_{k} \), \( \alpha_{k} \)
are non-negative integers, \( \{j_{1},\ldots,j_{r}\} \subset \{1,\ldots,n\} \),
\(
K_{i} \) are compact sets satisfying
\( K_{i} \) lies in the interiors of \( K_{i+1} \) and
\( \cup_{i=1}^{\infty} K_{i} = \Omega \). As \( |\alpha|=0 \),
\( \partial^{|\alpha|} f/ \partial
x_{j_{1}}^{\alpha_{1}}
\cdots  \partial x_{j_{r}}^{\alpha_{r}}=f \).
Endowed with this metric topology, \( C^{\infty}(\Omega) \) is a locally convex
space with a complete translation invariant metric, i.e., a Fr\'echet space, and
\( [D(\Omega),d_{C^{\infty}(\Omega)}] \) is a translation invariant metric vector space.

\( B_{d}[D(\Omega),C] \) is a Banach space and
a unital Banach algebra by Corollary 11 and Corollary 12. Let \( L_{p}(\Omega),
1 \leq p < \infty \), be the spaces of complex-valued functions \( x
\) defined on \( \Omega \) satisfying that
\( |x|^{p} \) is integrable with respect to the Lebesgue measure. Let
\( m \) of which support is a subset of \( K_{N} \) be a measurable complex function
defined on \( \Omega \). Then, if \( m
\) is a Lebesgue integrable complex function, then
\( \Lambda_{m} \), the corresponding linear functional, i.e., a distribution
with respect to another topology on \( D(\Omega) \), falls in the space \(  B_{d}[D(\Omega),C] \)
owing to
\begin{eqnarray*}
&& \left|\left|\Lambda_{m}\right|\right|_{B_{d}[D(\Omega),C]}\\ &=&
\sup_{f \neq
0,f \in D(\Omega)}\frac{\left|\Lambda_{m}(f)\right|}
{d_{C^{\infty}(\Omega)}(f,0)}\\ &\leq& \sup_{f \neq 0,f \in
D(\Omega)}\frac{p_{N}(f)\left|\left|m\right|\right|_{L_{1}(\Omega)}}
{d_{C^{\infty}(\Omega)}(f,0)}\\ &\leq&
\left|\left|m\right|\right|_{L_{1}(\Omega)}.
\end{eqnarray*}
 If \( m \in L_{\infty}(\Omega) \), the corresponding linear functional \( \Lambda_{m} \)
also falls in the space \( B_{d}[D(\Omega),C] \), where \( L_{\infty}(\Omega)
\) is the space of all essentially bounded functions on \( \Omega \). If
\( m \in L_{p}(\Omega) \), \( 1 < p < \infty \), and the Lebesgue
measure on the set \( \{x: m(x) \leq 1,x \in \Omega\} \) is finite, the corresponding linear
functional \( \Lambda_{m} \) falls in the space \(B_{d}[D(\Omega),C] \).

In addition to the linear functionals corresponding to the "ordinary" functions,
consider
\( \delta \), the Dirac delta function as a linear functional on
\( D(\Omega) \). Then  \( \delta_{c} \), \( c \in K_{N} \),
\begin{eqnarray*}
&& d(\delta_{c},0) \\ &=& \sup_{f \neq 0,f \in
D(\Omega)}\frac{\left|f(c)\right|} {d_{C^{\infty}(\Omega)}(f,0)} \\ &\leq& 1,
\end{eqnarray*}
i.e., \( \delta_{c} \in B_{d}[D(\Omega),C] \), where \( \delta_{c}(f)=f(c) \) for \( f
\in D(\Omega) \). Analogously, as \( 0 \in K_{N} \),
\begin{eqnarray*}
&& d(\delta^{(k)},0) \\ &=& \max\left\{\sup_{f \neq 0,f \in
D(\Omega)}\frac{\left|\delta^{(k)}(f)\right|}
{d_{C^{\infty}(\Omega)}(f,0)},\left|\delta^{(k)}(0)\right|\right\} \\ &=& \sup_{f
\neq 0,f
\in D(\Omega)}\frac{\left|f^{(k)}(0)\right|} {d_{C^{\infty}(\Omega)}(f,0)} \\ &\leq& 1,
\end{eqnarray*}
i.e., \( \delta^{(k)} \in B_{d}[D(\Omega),C] \), where \( k=(k_{1},k_{2},\ldots,k_{r}) \), \(
\delta^{(k)}(f)=(-1)^{k}f^{(k)}(0)
\)  for \( f
\in D(\Omega) \),
\( f^{(k)}=\partial^{|k|} f/\partial
x_{j_{1}}^{k_{1}}
\cdots  \partial x_{j_{r}}^{k_{r}} \),
and \( |k|
\leq N
\). \\

If the metric \( d_{D(\Omega) } \) imposed on \( D(\Omega) \) (not on \( C^{\infty}(\Omega) \)) has
the same form as \( d_{C^{\infty}(\Omega) } \) with \( p_{i} \) modified to
\begin{eqnarray*}
p_{i}(f)=\max\left\{\left|\left(\frac{\partial^{|\alpha|} f} {\partial
x_{j_{1}}^{\alpha_{1}}
\cdots  \partial x_{j_{r}}^{\alpha_{r}}}\right)(x)\right|:
x=(x_{1},\ldots,x_{n})^{t} \in \Omega,
\left|\alpha\right| \leq i\right\},
\end{eqnarray*}
then \( \delta_{c} \) with \( c \in \Omega \) and
\( \delta^{(k)} \) with \( 0 \in \Omega \), \( |k| \leq N \),  fall in the space \(B_{d}[D(\Omega),C] \) with respect to
this metric. Furthermore, without the assumption that the support of \( m \) is a subset of \( K_{N} \),
\( \Lambda_{m} \in B_{d}[D(\Omega),C] \) still holds for \( m \in L_{1}(\Omega) \)
and for \( m \in L_{p}(\Omega) \) with the Lebesgue
measure on the set \( \{x: m(x) \leq 1,x \in \Omega\} \) being finite, \( 1 <  p < \infty \). As the volume of
\( \Omega \) is finite and \( m \in L_{\infty}(\Omega) \), the corresponding linear functional
\( \Lambda_{m} \in B_{d}[D(\Omega),C] \) with respect to the metric  \( d_{D(\Omega) } \).

The above linear functionals falling in the space \( B_{d}[D(\Omega),C] \)
are also continuous by Theorem 9. In addition, as \( \Lambda_{m} \), \( \delta_{c} \),
and \( \delta^{(k)} \) fall in \( B_{d}[D(\Omega),C] \), the square of these mappings, i.e.,
\( \Lambda_{m}*\Lambda_{m} \), \( \delta_{c}*\delta_{c} \),
and \( \delta^{(k)}*\delta^{(k)} \) with the multiplication operation given
in Corollary 12, also fall in \( B_{d}[D(\Omega),C] \) and are not linear.

\end{example}

\begin{example}

Let \( \Omega \subset R \) be a nonempty open
set. Then, \( B_{d}[D(\Omega),C]
\) (also see Example 1) is a Banach space and a unital Banach algebra by Corollary 11 and Corollary 12. \\
(a): Let \(  F:Dom(F) \rightarrow B_{d}[D(\Omega),C] \) be the operator
defined by
\begin{eqnarray*}
F(\Lambda)(\phi)=\Lambda(x\phi)
\end{eqnarray*}
for \( \Lambda \in Dom(F) \), where
\( Dom(F) \subset B_{d}[D(\Omega),C] \) consists of some functionals
\( \Lambda \) in \( B_{d}[D(\Omega),C] \) satisfying that \( F(\Lambda) \in
B_{d}[D(\Omega),C] \),
\( \phi \in D(\Omega) \), and
\( x \) is a real-valued function on \( \Omega \) defined by \( x(t)=t \) for
\( t \in \Omega \). As
\( \Lambda \) are the linear functionals corresponding to the
Dirac's delta function or the square integrable functions in \( L_{2}(\Omega)
\),
\( F \) is associated with the multiplication operator, i.e., being associated with the position operator in quantum mechanics. However,
\( F \) might not be associated with the multiplication operator as
\( \Lambda \) are some other nonlinear functionals in \( B_{d}[D(\Omega),C] \). If
\( Dom(F) \subset \{\delta_{c}: c \in K_{N}\} \), then \( F \in B\{Dom(F),B_{d}[D(\Omega),C]\} \) is
bounded with respect to the metric \( d_{C^{\infty}(\Omega) } \) owing to \( F(\delta_{c})=c\delta_{c} \) and thus
\begin{eqnarray*}
\left|\left|F\right|\right|_{B\left\{Dom(F),B_{d}\left[D(\Omega),C\right]\right\}}=
\sup_{\delta_{c} \in Dom(F)}\frac{\left|\left|F\left(\delta_{c}\right)\right|\right|_{B_{d}[D(\Omega),C]}}
{\left|\left|\delta_{c}\right|\right|_{B_{d}[D(\Omega),C]}}\leq \sup_{c \in K_{N}}\left|c\right|.
\end{eqnarray*}
However, as the metric on \( D(\Omega) \) is \( d_{D(\Omega)} \), \( \{\delta_{c}: c \in \Omega\} \subset Dom(F)\), and \( \Omega \) is unbounded,
then \( F \) is not bounded. \\
(b): Let \( F:Dom(F) \rightarrow B_{d}[D(\Omega),C]
\) defined by
\( F(\Lambda)(\phi)=(-1)^{\alpha}\Lambda(\phi^{(\alpha)}) \), where
\( Dom(F) \subset B_{d}[D(\Omega),C] \) consists of some functionals
\( \Lambda \) in \( B_{d}[D(\Omega),C] \) satisfying that \( F(\Lambda) \in
B_{d}[D(\Omega),C] \), \( \phi \in D(\Omega) \), and \( \alpha \)
is a positive integer. \( F \) is the differential operator on the "bounded"
linear functionals. Note that \( F \) is associated the momentum operator in
quantum mechanics as \( \alpha=1 \). Further, if \( \alpha=1 \),  \( F
\) is unbounded because
\begin{eqnarray*}
&&
\frac{\left|\left|F(\Lambda_{n})\right|\right|_{B_{d}[D(\Omega),C]}}
{\left|\left|\Lambda_{n}\right|\right|_{B_{d}[D(\Omega),C]}}
\geq 2n,
\end{eqnarray*}
where \( \Lambda_{n} \in Dom(F) \) are the linear functionals, i.e., the distributions on
\( D(\Omega) \) endowed with another topology, corresponding to the
functions
\( \lambda_{n}(t)=n(t-c), c \leq t \leq c+1/n \) and \( 0 \) otherwise, \( n \geq N_{0} \),
\( c \in \Omega \), and where \( N_{0} \) is some positive integer.

\end{example}

\begin{example}

(a): Let \( F: L_{1}(R^{n}) \rightarrow C_{0}(R^{n}) \) be the Fourier transform
defined by
\begin{eqnarray*}
\left[F(f)\right](t)=\hat{f}(t)=\left(2 \pi \right)^{-n/2}\int_{R^{n}}fe^{-it \cdot x}dx
\end{eqnarray*}
for \( f \in L_{1}(R^{n}) \), where \( C_{0}(R^{n}) \) endowed with the
supremum norm is the Banach space of all complex continuous functions on \(
R^{n} \) that vanish at infinity. Since
\( \|\hat{f}\|_{C_{0}(R^{n})} \leq \|f\|_{L_{1}(R^{n})} \)
(see \cite{rudin}, Theorem 7.5),
\begin{eqnarray*}
\left|\left|F\right|\right|_{B[L_{1}(R^{n}),C_{0}(R^{n})]}
=\sup_{f\neq 0,f \in L_{1}(R^{n})}
\frac{\left|\left|\hat{f}\right|\right|_{C_{0}(R^{n})}}
{\left|\left|f\right|\right|_{L_{1}(R^{n})}} \leq 1
\end{eqnarray*}
and hence \( F \in B[L_{1}(R^{n}),C_{0}(R^{n})] \).

As \( F \) is the Fourier-Plancherel transform, i.e.,
\( F: L_{2}(R^{n}) \rightarrow L_{2}(R^{n}) \), \( F \in
B[L_{2}(R^{n}),L_{2}(R^{n})] \) since it is a linear isometry. \\
(b): Let \( \mathcal{L}(R^{n}) \subset C^{\infty}(R^{n}) \) (see Example 1 and
Example 2) be the Fr\'echet space with the metric
\begin{eqnarray*}
d(f,g)=\max\left(\sum_{k=0}^{\infty}\frac{b^{-k}\left|\left|f-g\right|\right|_{k}}
{a+\left|\left|f-g\right|\right|_{k}},\left|\left|f-g\right|\right|_{N}\right)
\end{eqnarray*}
for \( f, g \in \mathcal{L}(R^{n}) \), where \( N \) is a positive integer, \( a >0 \), \( b > 1 \),
\begin{eqnarray*}
\left|\left|f\right|\right|_{k}=\sup_{\left|\alpha\right| \leq k}\sup_{x \in R^{n}}
\left(1+\left|\left|x\right|\right|_{R^{n}}^{2}\right)^{k}
\left|\left(\frac{\partial^{|\alpha|} f} {\partial
x_{j_{1}}^{\alpha_{1}}
\cdots  \partial x_{j_{r}}^{\alpha_{r}}}\right)(x)
\right|,
\end{eqnarray*}
and where \( ||\cdot||_{R^{n}} \) is the usual Euclidean norm on \( R^{n} \).
Note that \( ||f||_{k} \)  is finite for \( f \in \mathcal{L}(R^{n}) \).
\( B_{d}[\mathcal{L}(R^{n}),C] \) is a Banach space and
a unital Banach algebra by Corollary 11 and Corollary 12. The linear functionals corresponding to \( \delta_{c} \), \( c \in R^{n} \),
\( \delta^{(k)} \), \( |k| \leq N \), and the Lebesgue integrable complex function \( m \), i.e.,
the tempered distributions (see \cite{rudin}, Definition 7.11) with respect to some
topology on \( \mathcal{L} (R^{n}) \), fall in \( B_{d}[ \mathcal{L}(R^{n}),C] \).
If \( m \in L_{p}(R^{n}) \) with the Lebesgue
measure on the set \( \{x: m(x) \leq 1,x \in R^{n}\} \) being finite, \( 1 < p < \infty \), the corresponding linear
functional \( \Lambda_{m} \) falls in the space \( B_{d}[ \mathcal{L}(R^{n}),C] \).

As \( F: \mathcal{L}(R^{n}) \rightarrow \mathcal{L}(R^{n}) \) is the Fourier
transform,
\( F \) is bounded, i.e., mapping bounded sets in \( \mathcal{L}(R^{n}) \) into bounded sets in \( \mathcal{L}(R^{n}) \),
because \( F \) is a continuous linear mapping from
\( \mathcal{L}(R^{n}) \) into \( \mathcal{L}(R^{n}) \) (see \cite{rudin}, Theorem
7.4).

\end{example}

\end{document}